\documentclass[12pt]{article}
\usepackage{amssymb}

\newenvironment{dfn}{\bigskip \noindent \bf Definition \rm}{\bigskip}
\newenvironment{rem}{\bigskip \noindent \bf Remark \rm}{\bigskip}

\newenvironment{proof}{\bigskip \noindent \bf Proof: \rm}{\bigskip}
\newenvironment{enumdfn}{
\begin{enumerate}}{\end{enumerate}}
\newenvironment{enumthm}{
\begin{enumerate}}{\end{enumerate}}
\newtheorem{thm}{Theorem}
\newtheorem{lemma}[thm]{Lemma}
\newtheorem{cor}[thm]{Corollary}
\newcommand{\qed}{\begin{flushright} \vspace{-1pc} $\square$
\end{flushright}}
\newcommand{\sthree}{\mbox{$\mathbb S^3$}}
\newcommand{\zdeg}[1]{\mbox{\it z-degree$(\Gamma_{#1})$}}
\newfont{\knot}{flype scaled 1000}
\begin{document}

\begin{center}
{\Large\bf Tait's Flyping Conjecture \\[1ex]
for 4-Regular Graphs} \\[3ex]
\large J\"org Sawollek\footnote{Fachbereich Mathematik, Universit\"at
Dortmund, 44221 Dortmund, Germany \\ {\em E-mail:\/}
sawollek@math.uni-dortmund.de \\ {\em WWW:\/}
http://www.mathematik.uni-dortmund.de/lsv/sawollek} \\[1ex]
June 22, 1998 (revised: October 9, 2002)
\end{center}
\vspace{4ex}

\begin{abstract}
Tait's flyping conjecture, stating that two reduced, alternating, prime
link diagrams can be connected by a finite sequence of flypes, is
extended to reduced, alternating, prime diagrams of 4-regular graphs in
\sthree . The proof of this version of the flyping conjecture is based
on the fact that the equivalence classes with respect to ambient isotopy
and rigid vertex isotopy of graph embeddings are identical on the class
of diagrams considered. \\[1ex]
\\[1ex]
\end{abstract}

\section*{Introduction}
Very early in the history of knot theory attention has been paid to
alternating diagrams of knots and links. At the end of the 19th century
P.G. Tait \cite{tait} stated several famous conjectures on alternating
link diagrams that could not be verified for about a century. The
conjectures concerning minimal crossing numbers of reduced, alternating
link diagrams \cite[Theorems A, B]{mura} have been proved independently
by Thistlethwaite \cite{this}, Murasugi \cite{mura}, and Kauffman
\cite{kauf87}. Tait's flyping conjecture, claiming that two reduced,
alternating, prime diagrams of a given link can be connected by a finite
sequence of so-called {\em flypes\/} (see \cite[page 311]{hkw} for
Tait's original terminology), has been shown by Menasco and
Thistlethwaite \cite{flype}, and for a special case, namely, for
well-connected diagrams, also by Schrijver \cite{schr}.

The present article, as well as \cite{alt}, deals with generalizations
of Tait's conjectures to embeddings of 4-regular (topological) graphs
into 3-space. In \cite{alt} it has been shown that a reduced,
alternating graph diagram $D$ has minimal crossing number. Furthermore,
if $D$ is prime in addition, then a non-alternating diagram that is
equivalent to $D$ cannot have the same crossing number as $D$.

The purpose of this paper is to prove Tait's flyping conjecture for
reduced, alternating, prime diagrams of 4-regular graphs in the 3-sphere
\sthree . The result depends on a suitable definition of primality, see
section \ref{properties}. Its proof is based on the fact that the
equivalence classes with respect to rigid vertex isotopy and ambient
isotopy of graph diagrams are identical on the class of diagrams under
consideration.

Definitions of these two equivalence relations for graph diagrams are
given in section \ref{diagrams}. Then, in section \ref{tangleinv}, the
notion of {\em tangles\/} is introduced to derive invariants of graph
diagrams and to give, via transformation tangles, a description of a
sequence of Reidemeister V moves (see Fig. \ref{reidemeister}) applied
to a graph vertex. In section \ref{properties} certain properties of
graph diagrams are discussed, and two important theorems on reduced,
alternating diagrams are stated. After, in section \ref{mainsec}, Tait's
flyping conjecture has been shown to be an immediate consequence of the
fact that the two equivalence classes mentioned above coincide for
reduced, alternating, prime graph diagrams, the latter statement is
finally proved in section \ref{mainproof}.

\section{Diagrams of 4-Regular Graphs}
\label{diagrams}

Embedding graphs into \sthree\ extends, in a natural way, the classical
knot theoretical problem of embedding one or more disjoint copies of the
1-sphere $\mathbb S^1$ into \sthree\ where the resulting images are
called {\em knots\/} and {\em links}, respectively. Classical
terminology of knot theory can be found in \cite{bur} or \cite{rol}, see
\cite{kaufbook}, \cite{kawa}, \cite{lick}, \cite{mubook} for more recent
introductions to the field.

A {\em topological graph\/} is a 1-dimensional cell complex which is
related to an abstract graph in the obvious way. In the following,
always $4$-regular graphs, which are allowed to have multiple edges or
loops, are considered. Vertices of degree two may occur but are
neglected since they are uninteresting for a topological treatment.

\begin{figure}[htb]
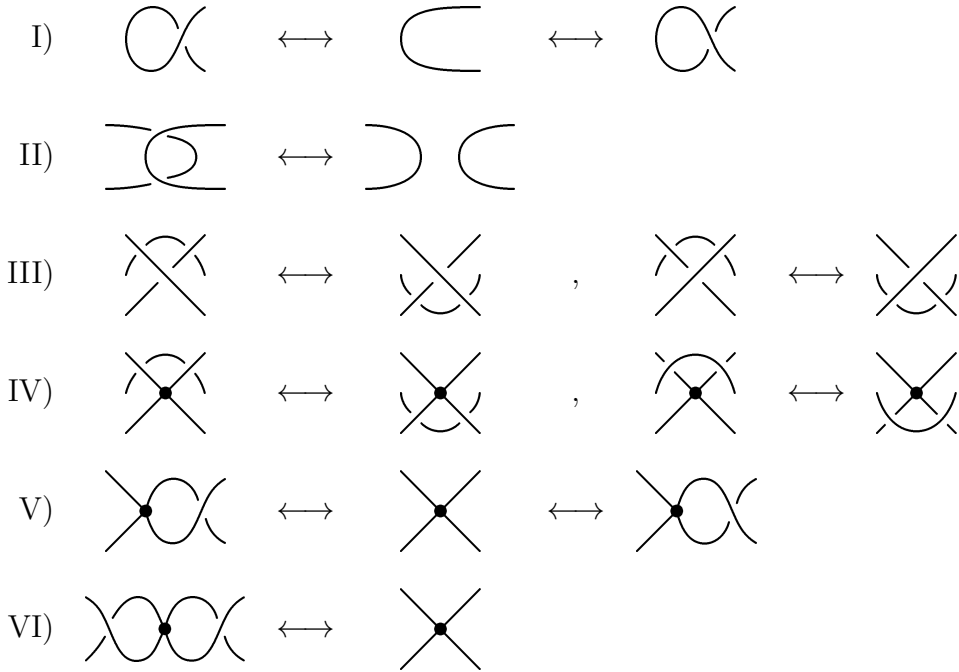

\begin{tabular}{rccccccc}
I) & {\knot a} & $\longleftrightarrow$ & {\knot b} &
$\longleftrightarrow$ & {\knot c} \\[2ex]
II) & {\knot d} & $\longleftrightarrow$ & {\knot e \hspace{0.5cm} f}
\\[2ex]
III) & {\knot g} & $\longleftrightarrow$ & {\knot h} & , & {\knot i} &
$\longleftrightarrow$ & {\knot j} \\[2ex]
IV) & {\knot k} & $\longleftrightarrow$ & {\knot l} & , & {\knot m} &
$\longleftrightarrow$ & {\knot n} \\[2ex]
V) & {\knot o} & $\longleftrightarrow$ & {\knot p} &
$\longleftrightarrow$ & {\knot q} \\[2ex]
VI) & {\knot r} & $\longleftrightarrow$ & {\knot p}
\end{tabular}
\vspace{2ex}

\caption{Reidemeister moves}
\label{reidemeister}
\end{figure}

If $G$ is a topological graph, then a {\em graph $\mathcal G$ in
\sthree} is the image of an embedding of $G$ into \sthree. Two graphs
${\mathcal G}_1$, ${\mathcal G}_2$ in \sthree\ are called {\em
equivalent with respect to ambient isotopy\/} or {\em ambient
isotopic\/} if there exists an orientation preserving autohomeomorphism
of \sthree\ which maps ${\mathcal G}_1$ onto ${\mathcal G}_2$.
Embeddings of topological graphs in \sthree\ can be examined via {\em
regular graph diagrams}, i.e., images under regular projections to an
appropriate sphere equipped with over-under information at double
points. Two graph diagrams $D$ and $D'$ are called {\em equivalent with
respect to ambient isotopy\/} or {\em ambient isotopic\/} if one can be
transformed into the other by a finite sequence of Reidemeister moves
I--V (see Fig. \ref{reidemeister}) combined with orientation preserving
homeomorphisms of the sphere to itself.
Two graphs in \sthree\ are ambient isotopic if and only if they have
diagrams that are ambient isotopic (see \cite{kauf89} or \cite{yet}).

Soon after the discovery of polynomial link invariants which fulfil
certain recurrence formulas, such as the Jones polynomial and its
generalizations, it had been tried to extend these invariants to graphs
in 3-space. Only quite recently such an invariant for arbitrary
topological graphs has been found by Yokota \cite{yoko} (see \cite{alt}
for a different approach in the case of 4-regular graphs). Yokota's
invariant manages the difficulty to be invariant under Reidemeister move
V which before had been the main obstacle for a full generalization of
combinatorial link invariants to knotted graphs.

Besides ambient isotopy there is a further equivalence relation for
graph diagrams which avoids Reidemeister move V and which will be
important for the purposes of this paper: two graph diagrams $D$ and
$D'$ are called {\em equivalent with respect to rigid vertex isotopy\/}
or {\em rigid vertex isotopic\/} if one can be transformed into the
other by a finite sequence of Reidemeister moves I--IV and VI (see Fig.
\ref{reidemeister}) combined with orientation preserving homeomorphisms
of the sphere to itself. Rigid vertex isotopy corresponds to an
equivalence relation on graphs in \sthree\ where a small neighbourhood
of each graph vertex is contained inside a disk, and only those
orientation preserving autohomeomorphisms of \sthree\ are considered
that respect these disks (see \cite{jomi}, \cite{kauf89}, \cite{yama}).
Observe that Reidemeister moves I--V imply Reidemeister move VI, thus
rigid vertex isotopic graph diagrams are ambient isotopic. See Fig.
\ref{counter} for two ambient isotopic diagrams which are not rigid
vertex isotopic.

For the sake of shortness, the phrase {\em graph diagram\/} will always
mean {\em (regular) diagram of a 4-regular graph in \sthree\/} in the
subsequent text, and throughout the article a {\em link\/} will be
considered as 4-regular graph in \sthree\ without vertices of degree
four. Of course, the equivalence classes with respect to rigid vertex
isotopy and ambient isotopy coincide for links.

\section{Tangles, and Invariants of Graph Diagrams}
\label{tangleinv}

In this section, the notion {\em rational tangle\/} is introduced for
two purposes: to define invariants of graph diagrams with respect to
ambient isotopy, and to describe the effect of a sequence of
Reidemeister moves applied to a graph vertex. The definitions and
notations used here are due to Conway \cite{con}.

\begin{figure}[htb]
\begin{picture}(12,2.5)
\thicklines
\put(0.5,0.8){\knot S}
\put(1.2,1.1){\line(0,1){0.45}}
\put(1.2,1.1){\line(1,0){0.25}}
\put(0.4,2){$a$}
\put(2.1,2){$b$}
\put(2.1,0.4){$c$}
\put(0.4,0.4){$d$}
\put(4.6,0.8){\knot T}
\put(5.3,1.1){\line(0,1){0.45}}
\put(5.3,1.1){\line(1,0){0.25}}
\put(7.45,1.1){\line(0,1){0.45}}
\put(7.45,1.1){\line(1,0){0.25}}
\put(5.3,2.1){$s$}
\put(7.5,2.1){$t$}
\put(6.1,0){$s+t$}
\put(8.9,0.8){\knot T}
\put(9.95,1.45){\line(-1,0){0.45}}
\put(9.95,1.45){\line(0,-1){0.25}}
\put(11.75,1.1){\line(0,1){0.45}}
\put(11.75,1.1){\line(1,0){0.25}}
\put(9.6,2.1){$s$}
\put(11.8,2.1){$t$}
\put(10,0){$s \cdot t = s t$}
\end{picture}
\caption{\label{tangle} A tangle, and tangle operations}
\vspace{1ex}

\end{figure}
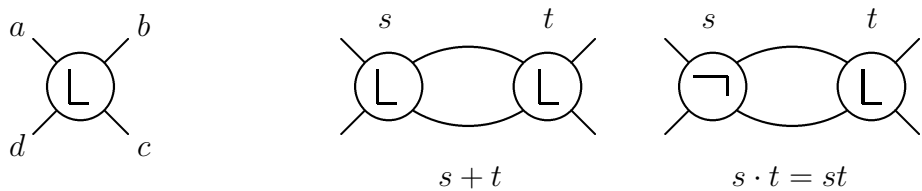

A {\em tangle\/} is a part of a link diagram in form of a disk with four
arcs emerging from it, see Fig. \ref{tangle} (left),
where the tangle's position is indicated by an L-shaped symbol and its
emerging arcs are labeled with letters $a$, $b$, $c$, $d$ in a clockwise
ordering (or simply one of them with "$a$"). Two tangles are said to be
{\em equivalent (with respect to ambient isotopy)\/} if one can be
transformed into the other by a finite sequence of Reidemeister moves of
type I--III and autohomeomorphisms of the disk which keep the boundary
fixed. In the following, a notational difference between a tangle and
the corresponding equivalence class will be avoided. Some basic tangles
are $0 = \raisebox{2pt}{\knot V}$, $\infty = \raisebox{2pt}{\knot W}$,
$1 = \mbox{\knot J}$, $\overline{1} = -1 = \mbox{\knot U}$.

For tangles $s$ and $t$, the operations + and $\cdot$ are defined as
depicted in Fig. \ref{tangle}. The tangles $0$, $n = 1 + \ldots + 1$,
and $\overline{n} = \overline{1} + \ldots + \overline{1}$ are called
{\em integer tangles}. If a tangle $t$ is of the form $t = a_1 \ldots
a_n$ with integer tangles $a_1, \ldots , a_n$ or if $t = \infty$, then
$t$ is called {\em rational tangle}. Let $\cal K$ denote the set of all
(equivalence classes of) rational tangles.

For a rational tangle $a_1 \ldots a_n$ the evaluation of the continued
fraction
\[ a_n + \frac{1}{a_{n-1} + \ldots + \frac{\textstyle 1}{\textstyle a_2 +
\frac{\textstyle 1}{\textstyle a_1}}} \]
gives a number in $\mathbb Q \cup \{ \infty \}$ where the obvious rules
for handling "$\infty$", such as $\frac{1}{\infty} = 0$, are applied, if
necessary, during the calculation. It is known that two rational tangles
are equivalent if and only if the values of the corresponding continued
fractions are identical (see \cite{con}, \cite{bur} for classical
proofs, and \cite{gold} for an elementary combinatorial proof).
Therefore, a rational tangle $r$ can be identified with this number,
thus let $r$ denote an element of $\cal K$ as well as the corresponding
value in $\mathbb Q \cup \{ \infty \}$, and let $|r|$ denote the
tangle's {\em crossing number}, i.e., the minimal number of crossings in
any diagram belonging to the equivalence class represented by $r$.
Furthermore, a rational tangle contained in ${\mathcal K} \setminus \{
0, \infty, 1, \overline{1} \}$ can be expressed in a uniquely determined
{\em normal form\/} $a_1 \ldots a_n$ such that $|a_1| \geq 2$, $a_2 \neq
0$, \ldots , $a_{n-1} \neq 0$, and all $a_i \geq 0$ or all $a_i \leq 0$.
The normal forms of the remaining four tangles are the obvious ones.
Tangles in normal form have minimal number of crossings.

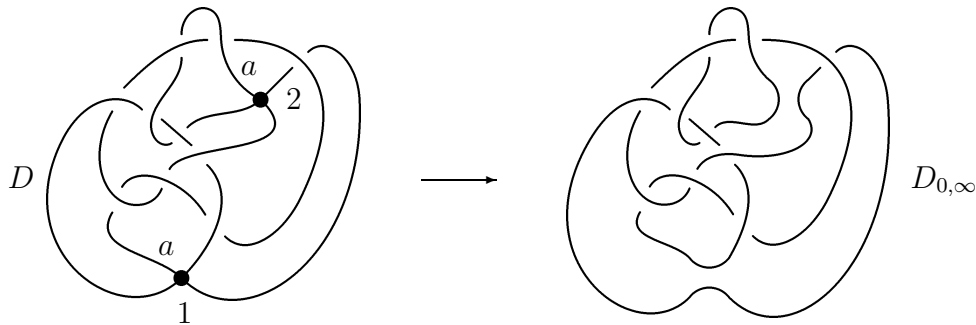
\begin{figure}[htb]
\begin{picture}(13.5,4.5)
\thinlines
\put(0.5,2){$D$}
\put(1.5,0.8){\knot X}
\put(3.6,3.5){$a$}
\put(4.2,3.05){2}
\put(2.5,1.1){$a$}
\put(2.75,0.2){1}
\put(6,2.1){\vector(1,0){1}}
\put(8.5,0.8){\knot Y}
\put(12.5,2){$D_{0, \infty}$}
\end{picture}
\caption{Cutting open vertices}
\label{replace}
\end{figure}

From a given graph diagram $D$ there can be obtained link diagrams by
substituting rational tangles for the graph vertices (see
\cite{four-reg}). To do this in a well-defined way, it is necessary
to give an ordering to the $k \geq 0$ graph vertices, i.e., a bijection
from the set $\{ 1 , \ldots , k \}$ to the set of graph vertices
contained in $D$ called {\em vertex-enumeration}, and an {\em
orientation\/} to each vertex, i.e., labeling an edge incident to the
vertex with the letter $a$. In the following, mainly the rational
tangles $0$, $\infty$, and $1$ will be needed. Substituting $0$ or
$\infty$ for graph vertices is done as depicted in Fig. \ref{replace},
that is to say, vertices are cut open in one of the two possible ways
determined by the vertex-orientation. Substituting the tangle $1$
corresponds to replacing a graph vertex with a crossing {\knot J} with
respect to the vertex-orientation given.

Replacing vertices of a graph diagram with all rational tangles, or, to
be precise, with a representative of each equivalence class (e.g., all
rational tangles in normal form),  gives an invariant of diagrams with
respect to ambient isotopy that consists of infinitely many link
diagrams (see \cite{four-reg}).

Getting invariants with respect to rigid vertex isotopy is much easier:
just substiute one or more arbitrary (but fixed) tangles for each vertex
and get rid of the ambiguity arising from different vertex-orientations
by considering all choices of such orientations. It is readily checked
that this gives well-defined invariants of graph diagrams under
Reidemeister moves I--IV and VI. For example, the unordered tuple
$(D_{0,0}, D_{0, \infty}, D_{\infty, 0}, D_{\infty, \infty})$ defines a
rigid vertex invariant of the graph belonging to the diagram $D$
depicted in Fig. \ref{replace}. Observe that the invariant $C(G)$
defined in \cite{kauf89} and denoted by ${\cal C}(G)$ in \cite{kauvo} is
a special version of this type of invariants, induced by the tangles $0$
and $1$, where sets are used instead of -- more informative -- unordered
tuples.

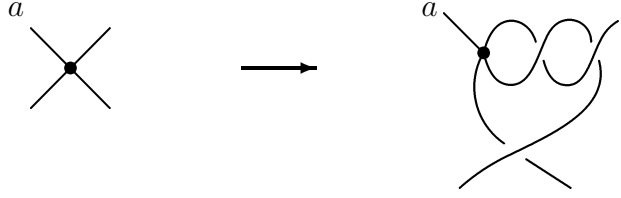
\begin{figure}[htb]
\begin{picture}(12,2.5)
\thicklines
\put(3.5,1.3){\knot p}
\put(3.2,2.1){$a$}
\put(6.3,1.4){\vector(1,0){1}}
\put(8.7,2.1){$a$}
\put(9,1.5){\knot F}
\end{picture}
\caption{Applying a sequence of Reidemeister V moves}
\label{sequence}
\end{figure}

Another point of view in considering ambient isotopic graph diagrams is
to observe that the Reidemeister moves of type V, applied to a vertex
during a transformation of one graph diagram into another, can be
collected to a rational tangle. For example, the transformation depicted
in Fig. \ref{sequence}
can be described by the tangle $t = 2 \, 1 \, 0 \in {\mathcal K}$.

\begin{dfn}
Let $D$, $D'$ be ambient isotopic graph diagrams with $k \geq 1$
vertices and given vertex-enumerations and -orientations. If $D_{0,
\ldots , 0} = D'_{t_1, \ldots , t_k}$ with tangles $t_1, \ldots, t_k$
then $t_j$, for $j \in \{ 1, \ldots , k \}$, is called {\em $j$--th
transformation tangle\/} of the transformation from $D$ into $D'$. If
$r_1, \ldots , r_k$ are rational tangles then let $r_j * t_j$ denote the
tangle into which $r_j$ is transformed when replacing the $j$--th vertex
of $D$ with $r_j$, that is to say, $D_{r_1, \ldots , r_k} = D'_{r_1 *
t_1, \ldots , r_k * t_k}$.
\end{dfn}

\begin{figure}[htb]
\begin{picture}(12,4.5)
\put(2.8,2.5){\begin{tabular}{cccccccc}
{\knot o} & $\longleftrightarrow$ & {\knot p} & $\longleftrightarrow$
& {\knot q} \\[2ex]
{\knot Q} & $\longleftrightarrow$ & {\knot p} & $\longleftrightarrow$ &
{\knot R}
\end{tabular}}
\put(5.9,4){$a$}
\put(2.7,4){$a$}
\put(8.6,4){$a$}
\put(5.9,1.85){$a$}
\put(3,1.85){$a$}
\put(8.8,1.85){$a$}
\end{picture}
\caption{Admissible Reidemeister moves of type V}
\label{admissible}
\end{figure}
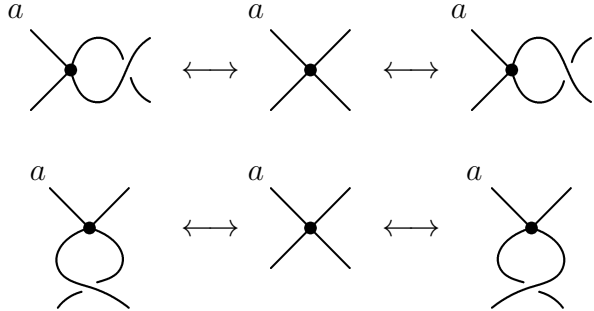

\begin{figure}[htb]
\begin{picture}(12,2)
\thicklines
\put(0.5,1.4){\knot C}
\put(2,1.5){\vector(1,0){1}}
\put(3.5,1.5){\knot K}
\put(7.5,1.4){\knot B}
\put(9,1.5){\vector(1,0){1}}
\put(10.5,1.5){\knot L}
\end{picture}
\caption{Replacing a graph vertex with 0 and $\infty$}
\label{trivialreplacement}
\end{figure}
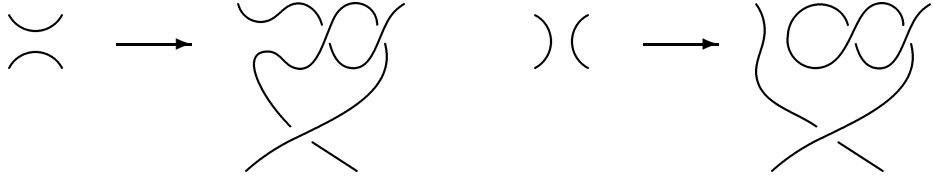

Considering Reidemeister move VI, it may be assumed, without loss of
generality, that only the {\em admissible\/} Reidemeister moves of type
V depicted in Fig. \ref{admissible} have to be applied during a
transformation.
Furthermore, if the orientation of a vertex is chosen appropriately, the
corresponding transformation tangle $t$ can be written as $t = b_1
\ldots b_s \neq \infty$ with integer tangles $b_1, \ldots , b_s$ such
that $b_1 \ldots b_s$ is in normal form. A rational tangle $r = a_1
\ldots a_l$ is transformed by the transformation tangle $t$ into the
tangle  $r*t = a_1 \ldots a_{l-1} (a_l+b_1) b_2 \ldots b_s$. Especially,
the tangles $0$ and $\infty $ are transformed into $t$ and into a tangle
equivalent to $b_2 \ldots b_s$, respectively, see Fig.
\ref{trivialreplacement}.

\begin{rem}
It is not difficult to see that the inverse of a transformation
described by a tangle $t = b_1 \ldots b_s$ is given by
$t' = \overline{b_s \ldots b_1}$, i.e., if $D_0 = D'_t$ then $D'_0 =
D_{t'}$ for ambient isotopic graph diagrams $D$ and $D'$ in which
corresponding vertices have been replaced.
\end{rem}

\section{Properties of Graph Diagrams}
\label{properties}

A link diagram is said to be {\em alternating\/} if over- and
undercrossings are alternating with each other while walking along any
link component in the diagram. A link diagram $D$ is called {\em
reduced\/} if it contains no {\em isthmus}, i.e., a crossing $p$ such
that $D \setminus \{ p \}$ has more components than $D$. A connected
link diagram is said to be {\em prime\/} if it cannot be written as
connected sum of two link diagrams both of which having at least one
crossing.

\begin{dfn}
A graph diagram $D$ is said to be {\em alternating/reduced/prime\/} if,
corresponding to an arbitrarily chosen vertex-enumeration and
-orientation, the link diagrams $D_{i_1, \ldots , i_k}$ are
alternating/reduced/prime for every choice of $i_1, \ldots , i_k \in \{
0, \infty \}$.
\end{dfn}

An example of a graph diagram that is alternating, reduced, and prime is
depicted in Fig. \ref{replace}.

\begin{rem}
\begin{enumerate}
\item The definition of primality for graph diagrams does not seem to be
the natural one, but it is the one that fits into the context (see the
counterexample contained in the remark at the end of this section).
\item In contrast to the case of link diagrams, a graph diagram that is
not reduced may be {\em irreducible}, i.e., the number of the diagram's
crossings is minimal. See \cite{alt} for examples.
\item A prime graph diagram which has more than one crossing is reduced.
\item The definition of "alternating" for graph diagrams has been
adapted to the definitions of "prime" and "reduced". Of course, a graph
diagram $D$ is alternating if and only if there is a choice of
vertex-orientations such that $D_{1, \ldots , 1}$ is alternating.
\end{enumerate}
\end{rem}

It is remarkable that, after introducing an appropriate definition of
primality (see \cite{alt}), for a 4-regular graph in \sthree\ which
possesses a reduced, alternating diagram the property to be prime can be
deduced from the corresponding property of its diagram. A proof of this
fact for link diagrams is due to Menasco \cite{mena}, and it can easily
be extended to 4-regular graphs in \sthree, see \cite[Theorem 8]{alt}.
The following statement is an immediate consequence.

\begin{thm}
\label{graphmenasco}
Let $D$ and $D'$ be ambient isotopic graph diagrams that are alternating
and reduced. Then:
\begin{enumthm}
\item $D$ is connected if and only if $D'$ is connected.
\item $D$ is prime if and only if $D'$ is prime.
\end{enumthm}
\end{thm}
\qed

To prove Tait's flyping conjecture for graph diagrams, a generalization
of a Tait conjecture concerning minimal crossing numbers, cited in the
introduction, will be needed. A proof can be found in \cite[Theorem
9]{alt}.

\begin{thm}
\label{altmin}
Let $\mathcal G$ be a 4-regular graph in \sthree, and let $D$ be a
reduced, alternating, prime diagram of $\mathcal G$ with $n$ crossings.
Then there is no diagram of $\mathcal G$ having less than $n$ crossings,
and any non-alternating diagram of $\mathcal G$ has more than $n$
crossings.
\end{thm}
\qed

\section{Tait's Flyping Conjecture}
\label{mainsec}

The definition of a tangle as part of a link diagram can be extended to
graph diagrams in the obvious way, and appropriate equivalence relations
can be introduced likewise. In the following, always those generalized
tangles are considered.

\begin{dfn}
A {\em flype\/} is a local change in a graph diagram as depicted in Fig.
\ref{flype}.
\begin{figure}[htb]
\begin{picture}(13,2)
\thicklines
\put(2.3,0.5){\knot O}
\put(4.2,0.8){\line(0,1){0.45}}
\put(4.2,0.8){\line(1,0){0.25}}
\put(4.2,0){$t$}
\put(6.5,0.9){$\longleftrightarrow$}
\put(9,0.5){\knot P}
\put(9.4,0.8){\line(0,1){0.1}}
\put(9.4,0.95){\line(0,1){0.1}}
\put(9.4,1.1){\line(0,1){0.1}}
\put(9.4,1.25){\line(1,0){0.1}}
\put(9.55,1.25){\line(1,0){0.1}}
\put(9.4,0){$t_h$}
\end{picture}
\caption{A flype}
\label{flype}
\end{figure}
Using Conway's notation \cite{con}, this corresponds to a transformation
of the form $1 + t \leftrightarrow t_h + 1$ or $\overline{1} + t
\leftrightarrow t_h + \overline{1}$.
\end{dfn}

In \cite{flype} Tait's flyping conjecture is proved for diagrams of
4-regular rigid vertex graphs to obtain the validity of Tait's original
conjecture. This result is stated in the following theorem where, as
well as in the subsequent text, the phrase {\em sequence of flypes\/} is
an abbreviation for {\em sequence of flypes plus orientation preserving
autohomeomorphisms of the sphere}. It should be mentioned that the
notion of primality for graph diagrams used in \cite{flype} is different
from the one used in the present text. Indeed, a graph diagram that is
prime with respect to the definition given above is prime in the sense
of \cite{flype}, too, and thus the result from \cite{flype} can be
adopted here.

\begin{thm}
\label{rigidflype}
Let $D$ and $D'$ be rigid vertex isotopic graph diagrams that are
reduced, alternating, and prime. Then there exists a finite sequence of
flypes which transforms $D$ into $D'$.
\end{thm}
\qed

Considering ambient isotopy of graph diagrams, the desired extension of
Tait's flyping conjecture to 4-regular graphs can be deduced immediately
from the next theorem. The proof of the theorem will be given in section
\ref{mainproof}.

\begin{thm}
\label{main}
Let $D$ and $D'$ be graph diagrams that are reduced, alternating, and
prime. Then $D$ and $D'$ are equivalent with respect to ambient isotopy
if and only if they are equivalent with respect to rigid vertex isotopy.
\end{thm}
\medskip

\begin{cor}
\label{maincor}
Let $D$ and $D'$ be ambient isotopic graph diagrams that are reduced,
alternating, and prime. Then there exists a finite sequence of flypes
which transforms $D$ into $D'$.
\end{cor}
\qed

\begin{rem}
Theorem \ref{main} does not hold in general without assuming primality.
A counterexample is depicted in Fig. \ref{counter}: the two (alternating
and reduced) graph diagrams obviously are equivalent with respect to
ambient isotopy, but they are not rigid vertex isotopic because cutting
open vertices in the two possible ways gives diagrams of 2- and
3-component links, respectively, where the diagrams belonging to the
3-component link are equivalent, and the diagrams belonging to the
2-component link correspond to different mirror images of a chiral link.
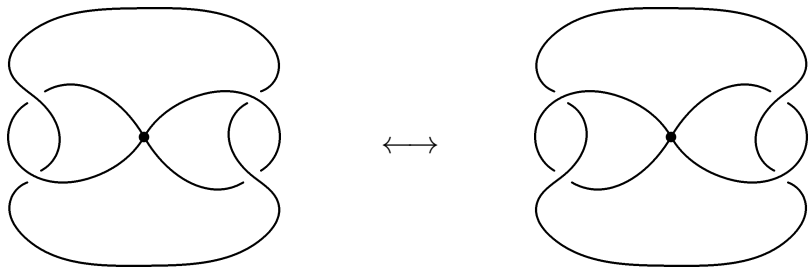
\begin{figure}[htb]
\begin{picture}(13,3.5)
\put(1.5,0){\knot H}
\put(6.5,1.5){$\longleftrightarrow$}
\put(8.5,0){\knot I}
\end{picture}
\caption{Ambient isotopic diagrams that are not rigid vertex isotopic}
\label{counter}
\end{figure}
\end{rem}

\section{Proof of Theorem \ref{main}}
\label{mainproof}

Theorem \ref{main} is proved by induction on the number of graph
vertices. The main ingredient for the induction step comes from the fact
that transformation tangles which describe a transformation between
reduced, alternating, prime diagrams always are trivial, and thus
replacing a vertex in both diagrams with the same tangle, for
corresponding vertices and with respect to appropriately chosen
vertex-orientations, gives ambient isotopic graph diagrams with one
vertex less. As a main tool for showing this claim, properties of the
{\em Kauffman polynomial\/} \cite{kauf90} of link diagrams and its
relations to a diagram's crossing number are used.

\begin{dfn}
The {\em Kauffman polynomial\/} $\Gamma_D(a,z) \in \mathbb Z [a^{\pm 1},
z^{\pm 1}]$ of a link diagram $D$ is defined by the following
properties.
\begin{enumdfn}
\item $\Gamma_D(a,z)=1$ if $D$ is a simple closed curve
\item $\Gamma_{D_{\mbox{\knot s v}}} = a \Gamma_{D_{\mbox{\knot u}}}$
and $\Gamma_{D_{\mbox{\knot t v}}} = a^{-1} \Gamma_{D_{\mbox{\knot u}}}$
\item $\Gamma_{D_{\mbox{\knot s}}} + \Gamma_{D_{\mbox{\knot t}}} \; = \;
z ( \Gamma_{D_{\mbox{\knot w}}} + \Gamma_{D_{\mbox{\knot x}}})$
\end{enumdfn}
The highest exponent in the variable $z$ is called {\em z-degree\/} of
$\Gamma_D$.
\end{dfn}

In the following, an {\em $n$-bridge\/} $b$ in a link diagram $D$ is an
arc of $D$ that contains only overcrossings or only undercrossings, and
the number of these crossings is $n$, the {\em length\/} of $b$.

The $z$-degree of the Kauffman polynomial heavily depends on the length
of the longest bridge contained in a diagram. A precise formulation of
this fact is given in the next theorem, for a proof see \cite[Theorems
4, 5]{this}.

\begin{thm}
\label{elementary}
Let $D$ be a link diagram with $n \geq 1$ crossings, and let $b \geq 1$
be the length of an arbitrary bridge contained in $D$. Then:
\begin{enumthm}
\item $\zdeg{D} \leq n - b \leq n - 1$
\item $\zdeg{D} = n - 1$ if and only if $D$ is reduced, alternating,
prime
\end{enumthm}
\end{thm}
\qed

Now two technical lemmas on Kauffman polynomials are stated which are
needed to deduce the crucial Lemma \ref{transform}. For a proof of the
first lemma see \cite[Lemmas 8, 10]{four-reg}.

\begin{lemma}
\label{tanglesum}
Let $D$ be a reduced, alternating, prime graph diagram with \mbox{$k
\geq 1$} vertices and $n$ crossings, and let $r$ be a rational tangle in
normal form. Furthermore, let $i_1, \ldots, i_{k-1} \in \{ 0, \infty
\}$. Then, corresponding to an arbitrarily chosen vertex-enumeration and
-orientation, the following holds for the link diagram $D' = D_{i_1,
\ldots, i_j, r, i_{j+1}, \ldots, i_{k-1}}$ with $j \in \{ 1, \ldots, k
\}$.
\[ \zdeg{D'} \; = \; \left\{
\begin{array}{ll}
n + |r| - 1 & \mbox{if $D'$ is alternating} \\[1ex]
n + |r| - 2 & \mbox{if $D'$ is not alternating and } |r| \neq 1
\end{array} \right. \]
\end{lemma}
\qed

\begin{dfn}
Let $D$ be an alternating graph diagram with $k \geq 1$ vertices and $n$
crossings, supplied with an arbitrarily chosen vertex-enumeration and a
vertex-orientation such that $D_{\overline{1}, \ldots, \overline{1}}$ is
alternating. Then $D$ is called {\em degree-reducing\/} if
$\zdeg{D_{\varepsilon_1, \ldots, \varepsilon_k}} \leq n-2$ for every
choice of $\varepsilon_1, \ldots, \varepsilon_k \in \{ 0, \infty, 1 \}$
with $\varepsilon_j = 1$ for at least one index $j \in \{ 1, \ldots, k
\}$.
\end{dfn}

\begin{lemma}
\label{ratins}
Let $D$ be a reduced, alternating, prime, degree-reducing graph diagram
with $k \geq 1$ vertices and $n$ crossings, supplied with an arbitrarily
chosen vertex-enumeration and a vertex-orientation such that
$D_{\overline{1}, \ldots, \overline{1}}$ is alternating. If $r_1, \ldots
r_k$ are rational tangles in normal form having at least two crossings
each, then
\[ \zdeg{D_{r_1, \ldots, r_k}} \; = \; n + |r_1| + \ldots + |r_k| - t_+
- 1 \]
holds, where $t_+$ denotes the number of indices $j$ with $r_j > 0$.
\end{lemma}

\begin{proof} As a consequence of Lemma \ref{tanglesum}, replacing a
vertex of $D$ with a negative tangle $r$ yields a graph diagram
$D'$ with $k-1$ vertices and $n + |r|$ crossings that is reduced,
alternating, prime ($\zdeg{D'_{\varepsilon_1, \ldots,
\varepsilon_{k-1}}} = n + |r| - 1$ with $\varepsilon_i \in \{ 0,
\infty\}$ can only be fulfilled if $D'$ possesses all three properties).
Thus it may be assumed that $t_+ = k$. In the following, it is shown
that
\[ \zdeg{D_{r_1, \ldots, r_k}} \left\{
\begin{array}{ll}
= n + |r_1| + \ldots + |r_k| - k - 1 & \mbox{if $r_j \neq 1$ for all
$j$} \\[1ex]
\leq n + |r_1| + \ldots + |r_k| - k - 2 & \mbox{otherwise}
\end{array} \right. \]
holds for positive rational tangles $r_1, \ldots, r_k$ having at least
one crossing each. Since Lemma \ref{tanglesum} and Theorem
\ref{elementary} (if $r_k = 1$) immediately give the result for the case
$k=1$, let $k \geq 2$ for the rest of the proof.

Let $l$ denote the number of indices $j$ with $r_j=1$. The proof is done
by induction on $k-l \geq 0$. If $k-l = 0$ then $r_1 = \ldots = r_k =
1$, and therefore $\zdeg{D_{r_1, \ldots, r_k}} \leq n-2$ since $D$ is
degree-reducing. Thus let $k-l \geq 1$. Then $r_j \neq 1$ for at least
one index $j \in \{ 1, \ldots , k \}$, and without loss of generality
let $r_k = a_1 \ldots a_s \neq 1$.

At first consider the case that $r_k = 2$. Then the recurrence formula
for the Kauffman polynomial gives:
\[ \Gamma_{D_{r_1, \ldots , r_{k-1}, 2}} \, = \, \Gamma_{D_{r_1, \ldots ,
r_{k-1}, 0}} + z(a \Gamma_{D_{r_1, \ldots , r_{k-1}, \infty}} +
\Gamma_{D_{r_1, \ldots , r_{k-1}, 1}}) \]
Applying the induction hypothesis immediately yields the desired
inequality if $l \geq 1$ and the desired equality if $l=0$.

If $r_k = a_1$ is integral then the claimed result can be verified
inductively by considering the corresponding recurrence formula, and a
further induction on $s$ completes the proof.
\qed
\end{proof}

An important class of degree-reducing diagrams is defined next.

\begin{dfn}
An alternating graph diagram $D$ with $k \geq 1$ vertices is called {\em
vertex-separating\/} if there exist disjoint tangles $t_1, \ldots, t_k$
in $D$ such that each tangle $t_i$ contains exactly one graph vertex and
replacing this vertex with a crossing {\knot J}, corresponding to an
appropriately chosen vertex-orientation, yields a 3-bridge inside $t_i$.
\end{dfn}

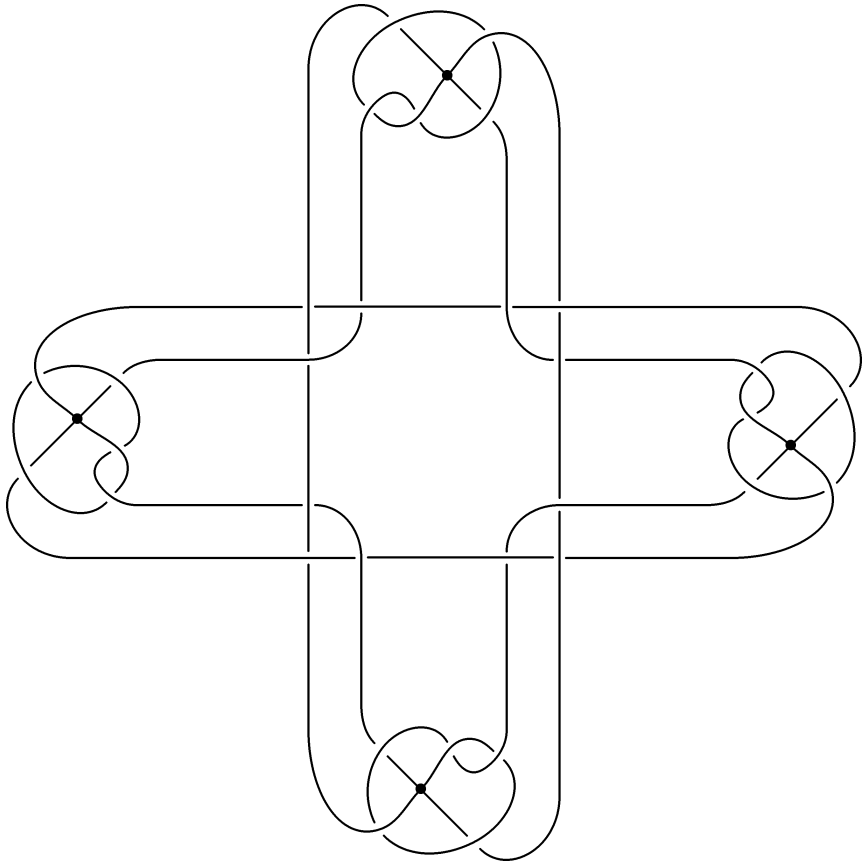
\begin{figure}[htb]
\begin{picture}(12,11.5)
\put(5.3,0.5){\knot M}
\put(5.3,0.5){\knot N}
\end{picture}
\caption{A reduced, alternating, prime, vertex-separating graph diagram}
\label{vertexsep}
\end{figure}

An example of an alternating graph diagram that is vertex-separating is
depicted in Fig. \ref{vertexsep}.

\begin{lemma}
A vertex-separating graph diagram is degree-reducing.
\end{lemma}

\begin{proof} Considering the {\em numerator formula for the Kauffman
polynomial\/} that has been deduced in \cite{four-reg}, p. 733, it can
easily be shown via induction on $l$ that
\[ \zdeg{D_{\varepsilon_1, \ldots, \varepsilon_k}} \; \leq \; n - l - 1
\; \leq \; n - 2 \]
holds for a vertex-separating graph diagram $D$ with $k$ vertices and
$n$ crossings where $\varepsilon_1, \ldots, \varepsilon_k \in \{ 0,
\infty, 1 \}$ have been replaced for the vertices and $l \geq 1$ denotes
the number of indices $j$ with $\varepsilon_j = 1$.
\qed
\end{proof}

\begin{lemma}
\label{transform}
Let $D$ and $D'$ be ambient isotopic graph diagrams with $k \geq 1$
vertices, and let $t_1, \ldots, t_k$ denote the transformation tangles
belonging to a sequence of Reidemeister moves which realizes the
equivalence between $D$ and $D'$, corresponding to chosen
vertex-enumerations and -orientations of the diagrams. If $D$ and $D'$
both are reduced, alternating, and prime then $|t_1| = \ldots = |t_k| =
0$.
\end{lemma}

\begin{proof} Obviously, a prime, reduced graph diagram with $k \geq 1$
vertices has at least two crossings, and, applying Theorem \ref{altmin},
it is clear that $D$ and $D'$ have identical crossing number $n \geq 2$.
Without loss of generality, let vertex-enumerations be chosen such that
the $j$-th vertex of $D$ is mapped to the $j$-th vertex of $D'$ when the
sequence of Reidemeister moves which transforms $D$ into $D'$ is
applied. For the sake of notational convenience, assume that
$D_{\overline{1}, \ldots, \overline{1}}$ and $D'_{\overline{1}, \ldots,
\overline{1}}$ are alternating.

Suppose that some of the transformation tangles, which may assumed to be
in normal form, have more than one crossing. Define $\varepsilon_i :=
\infty$ if $|t_i| = 1$ and $\varepsilon_i := 0$ if $|t_i| \neq 1$. Then
on the one hand, $\zdeg{D_{\varepsilon_1, \ldots, \varepsilon_k}} = n-1$
since $D$ is reduced, alternating, prime. On the other hand,
$\zdeg{D'_{\varepsilon_1*t_1, \ldots, \varepsilon_k*t_k}} \geq n$, by
Lemma \ref{ratins}, since $\varepsilon_i*t_i = \infty$ if $|t_i| = 1$
(perform a Reidemeister move I) and $D'_{\varepsilon_1*t_1, \ldots,
\varepsilon_k*t_k}$ can be thought to arise from a vertex-separating
graph diagram, namely, the diagram in which the vertices $v_i$ with
$|t_i| \neq 1$ not have been replaced (if there is no 3-bridge after
replacing one of these vertices by a crossing then perform a flype
as depicted in Fig. \ref{vertexflype}).
\begin{figure}[htb]
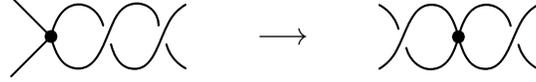

\[ \mbox{\knot 1} \qquad \longrightarrow \qquad \mbox{\knot 2} \]
\caption{Flyping a vertex}
\label{vertexflype}
\end{figure}
This gives a contradiction
because $D_{\varepsilon_1, \ldots, \varepsilon_k} =
D'_{\varepsilon_1*t_1, \ldots, \varepsilon_k*t_k}$. Thus $|t_i| \leq 1$
for all indices $i$.

Now suppose there is an index $j$ with $|t_j| = 1$. Define
$\varepsilon_j := 0$ and $\varepsilon_i := \infty$ for $i \neq j$. Then
$\zdeg{D_{\varepsilon_1, \ldots, \varepsilon_k}} = n-1$ as in the
previous case, but it is $\varepsilon_i * t_i = \infty$ for $i \neq j$
and therefore either $\zdeg{D'_{\varepsilon_1*t_1, \ldots,
\varepsilon_k*t_k}} \leq n-2$ if $t_j = 1$ because the diagram contains
a 3-bridge (otherwise there would be a contradiction to reducedness and
primality of the diagram), or $\zdeg{D'_{\varepsilon_1*t_1, \ldots,
\varepsilon_k*t_k}} = n$ if $t_j = \overline{1}$ because the diagram is
reduced, alternating, and prime. Again a contradiction in both cases,
thus $|t_i| = 0$ for all indices $i$.
\qed
\end{proof}

Lemma \ref{transform} shows that transformation tangles belonging to two
ambient isotopic graph diagrams always are trivial if both diagrams are
reduced, alternating, and prime. Thus replacing vertices of such
diagrams with tangles gives an ambient isotopy invariant of graph
diagrams up to a choice of vertex-orientations.
Especially:

\begin{cor}
\label{trafocor}
Let $D$ and $D'$ be ambient isotopic graph diagrams with $k \geq 1$
vertices that are reduced, alternating, and prime. Let $D_1$, $D'_1$
denote diagrams that arise from substituting the tangle $1$ for
corresponding vertices of $D$ and $D'$, respectively, and let
vertex-orientations be chosen such that $D_1$ and $D_1'$ both are
alternating. Then $D_1$ and $D'_1$ are equivalent with respect to
ambient isotopy.
\qed
\end{cor}

\begin{rem}
Observe that it is not clear, at this stage, that the diagrams $D_1$ and
$D'_1$ of Corollary \ref{trafocor} are rigid vertex isotopic since the
sequence of Reidemeister moves of type V applied to a vertex leads to a
transformation tangle which is equivalent to a trivial tangle but not
necessarily equal to it.
\end{rem}

\begin{lemma}
\label{rot}
The graph diagrams that are depicted in Fig. \ref{rotate}
\begin{figure}[htb]
\begin{picture}(8,2)
\thicklines
\put(1.8,0.5){\knot G}
\put(4.85,0.8){\line(0,1){0.45}}
\put(4.85,0.8){\line(1,0){0.25}}
\put(6.5,0.9){$\longleftrightarrow$}
\put(8,0.5){\knot G}
\put(11.4,1.15){\line(0,-1){0.25}}
\put(10.95,0.9){\line(1,0){0.45}}
\end{picture}
\caption{Rigid vertex isotopic diagrams}
\label{rotate}
\end{figure}
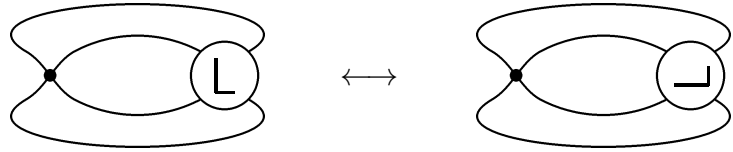
are rigid vertex isotopic.
\end{lemma}

\begin{figure}[htb]
\begin{picture}(8,9.5)
\thicklines
\put(-0.2,7.5){\knot G}
\put(2.85,7.8){\line(0,1){0.45}}
\put(2.85,7.8){\line(1,0){0.25}}
\put(3.8,7.9){$\longrightarrow$}
\put(4,8.1){\small 1}
\put(4.4,7.5){\knot y}
\put(7.45,7.8){\line(0,1){0.45}}
\put(7.45,7.8){\line(1,0){0.25}}
\put(8.4,7.9){$\longrightarrow$}
\put(8.6,8.1){\small 2}
\put(9,7.5){\knot A}
\put(12.05,7.8){\line(0,1){0.1}}
\put(12.05,7.95){\line(0,1){0.1}}
\put(12.05,8.1){\line(0,1){0.1}}
\put(12.05,8.25){\line(1,0){0.1}}
\put(12.2,8.25){\line(1,0){0.1}}
\put(1,3.9){$\longrightarrow$}
\put(1.2,4.1){\small 3}
\put(2.5,3.5){\knot z}
\put(5.55,3.8){\line(0,1){0.1}}
\put(5.55,3.95){\line(0,1){0.1}}
\put(5.55,4.1){\line(0,1){0.1}}
\put(5.55,4.25){\line(1,0){0.1}}
\put(5.7,4.25){\line(1,0){0.1}}
\put(4.7,4.9){\line(0,1){0.15}}
\put(4.7,5.15){\line(0,1){0.15}}
\put(4.7,5.4){\line(0,1){0.15}}
\put(4.7,5.65){\vector(0,1){0.4}}
\put(4.8,6.4){\circle*{0.1}}
\put(5,6.3){$\infty$}
\put(7.7,3.9){$\longrightarrow$}
\put(7.9,4.1){\small 4}
\put(9,4.1){\knot Z}
\put(12.05,4.4){\line(0,1){0.1}}
\put(12.05,4.55){\line(0,1){0.1}}
\put(12.05,4.7){\line(0,1){0.1}}
\put(12.05,4.85){\line(1,0){0.1}}
\put(12.2,4.85){\line(1,0){0.1}}
\put(1,0.9){$\longrightarrow$}
\put(1.2,1.1){\small 5}
\put(2.3,0.5){\knot G}
\put(5.65,1.15){\line(-1,0){0.1}}
\put(5.5,1.15){\line(-1,0){0.1}}
\put(5.35,1.15){\line(-1,0){0.1}}
\put(5.7,1.15){\line(0,-1){0.1}}
\put(5.7,1){\line(0,-1){0.1}}
\multiput(2.3,1)(0.4,0){7}{\line(1,0){0.2}}
\multiput(2.3,1)(0.4,0){7}{\line(1,0){0.2}}
\multiput(5.9,1)(0.4,0){2}{\line(1,0){0.2}}
\put(6.4,0.9){$\circlearrowleft$}
\put(7.3,0.9){$\longrightarrow$}
\put(8.3,0.5){\knot G}
\put(11.7,1.15){\line(0,-1){0.25}}
\put(11.25,0.9){\line(1,0){0.45}}
\end{picture}
\caption{1 = R. II; 2 = flype; 3 = R. VI; 4 = planar isotopy; 5 = R. II}
\label{picproof}
\end{figure}
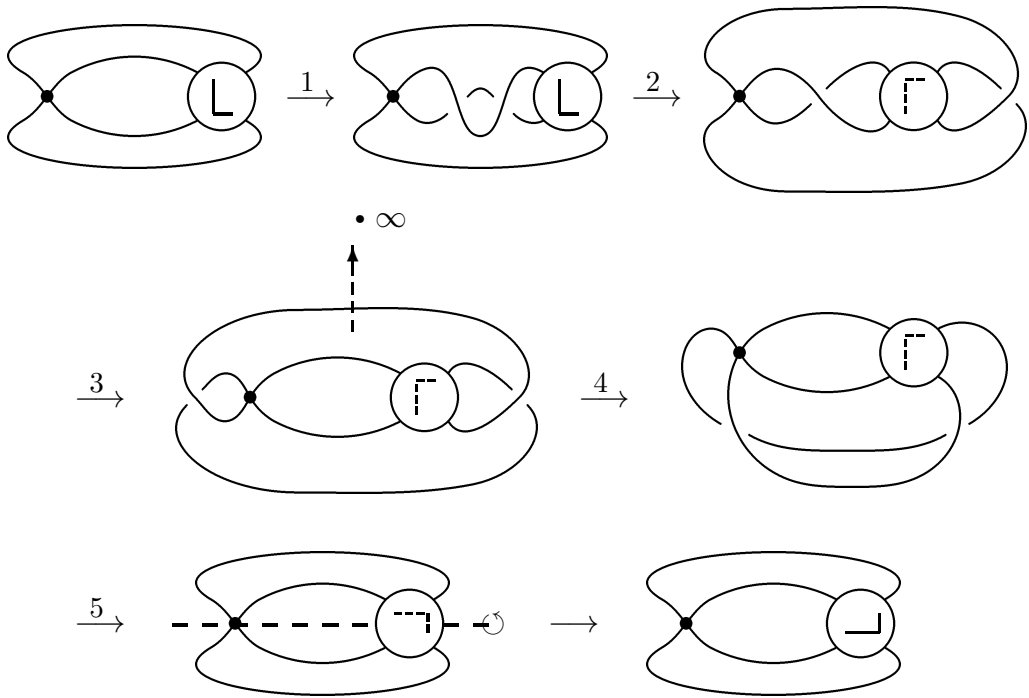

\begin{proof}
Perform the transformations depicted in Fig. \ref{picproof},
where "R." is an abbreviation for "Reidemeister", and observe that move
4 consists of pushing the diagram's upper arc upwards beyond "$\infty$".
\qed
\end{proof}

Now let $D$ and $D'$ be ambient isotopic diagrams which fulfil the
assumptions of Theorem \ref{main}. Of course, if the corresponding
graph has no vertices then the statement follows from Theorem
\ref{rigidflype}.

For the induction step, an arbitrary graph vertex of $D$ and the
corresponding vertex of $D'$ are replaced with a crossing each such that
the diagrams arising, which shall be denoted by $D_1$ and $D'_1$
respectively, are alternating. Then it follows from Corollary
\ref{trafocor} and the induction hypothesis that $D_1$ and $D'_1$ are
rigid vertex isotopic and thus can be connected by a finite sequence of
flypes, by Theorem \ref{rigidflype}.

Considering the defining Fig. \ref{flype}, there are three different
types of flype that can be applied to $D_1$: either the substituted
crossing lies outside the depicted part of the diagram, or it lies
inside the tangle $t$, or it is identical to the crossing next to $t$.
Call the latter one {\em essential\/} flype. Because of the obvious
one-to-one correspondence between the crossings before and after
applying a flype it is possible to keep track of the substituted
crossing during the whole sequence of flypes connecting $D_1$ and
$D'_1$. Obviously, if the sequence contains no essential flypes then it
gives rise to a sequence of flypes that can be applied to the diagram
$D$ and has $D'$ as final diagram, and the induction step is done.

\begin{figure}[htb]
\begin{picture}(13,2.5)
\thicklines
\put(0,0.5){\knot E}
\put(3.05,0.8){\line(0,1){0.45}}
\put(3.05,0.8){\line(1,0){0.25}}
\put(4.1,0.9){$\longrightarrow$}
\put(4.3,1.1){\small 1}
\put(5,0.5){\knot D}
\put(5.7,0.8){\line(0,1){0.1}}
\put(5.7,0.95){\line(0,1){0.1}}
\put(5.7,1.1){\line(0,1){0.1}}
\put(5.7,1.25){\line(1,0){0.1}}
\put(5.85,1.25){\line(1,0){0.1}}
\multiput(7.1,-0.25)(0,0.4){7}{\line(0,1){0.2}}
\put(6.95,2.3){$\circlearrowleft$}
\put(9,0.9){$\longrightarrow$}
\put(9.2,1.1){\small 2}
\put(9.7,0.5){\knot E}
\put(13,0.8){\line(0,1){0.45}}
\put(12.75,1.25){\line(1,0){0.25}}
\end{picture}
\caption{1 = flype; 2 = rotation around dashed axis}
\label{finalpic}
\end{figure}
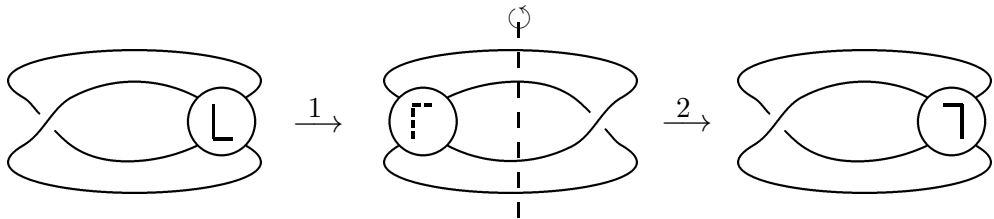

Now assume that the flyping sequence contains an essential flype. Since
$D$ is prime, cutting open the substituted crossing in either way must
yield a prime graph diagram, and therefore an essential flype can,
essentially, only be applied to $D_1$ as depicted in the first two
pictures of Fig. \ref{finalpic},
i.e., the whole diagram is involved in the flyping move. Applying Lemma
\ref{rot} twice shows that there arise rigid vertex isotopic diagrams
from re-inserting vertices for crossings in first and last diagram of
Fig. \ref{finalpic}. Therefore an essential flype applied to $D_1$ is
related to a rigid vertex isotopy of the diagram $D$, and an induction
on the number of essential flypes in the flyping sequence completes the
proof of Theorem \ref{main}.


\begin{thebibliography}{123}

\bibitem{bur} G. Burde and H. Zieschang, {\em Knots}, de Gruyter, Berlin
(1985).

\bibitem{con} J.H. Conway, An enumeration of knots and links, and some
of their algebraic properties, in: {\em Computational Problems in
Abstract Algebra\/} (ed. J. Leech), Pergamon Press, Oxford (1970),
329--358.

\bibitem{gold} J.R. Goldman and L.H. Kauffman, Rational Tangles, {\em
Adv. in Appl. Math.\/} {\bf 18} (1997), 300--332.

\bibitem{hkw} P. de la Harpe, M. Kervaire, and C. Weber, On the Jones
polynomial, {\em Enseign. Math.\/} {\bf 32} (1986), 271--335.

\bibitem{jomi} D. Jonish and K.C. Millett, Isotopy invariants of graphs,
{\em Trans. Amer. Math. Soc.\/} {\bf 327} (1991), 655--702.

\bibitem{kauf87} L.H. Kauffman, State models and the Jones polynomial,
{\em Topology\/} {\bf 26} (1987), 395--407.

\bibitem{kauf89} L.H. Kauffman, Invariants of graphs in three-space,
{\em Trans. Amer. Math. Soc.\/} {\bf 311} (1989), 697--710.

\bibitem{kauf90} L.H. Kauffman, An invariant of regular isotopy, {\em
Trans. Amer. Math. Soc.\/} {\bf 318} (1990), 417--471.

\bibitem{kaufbook} L.H. Kauffman, {\em Knots and Physics}, World
Scientific, Singapore (1991).

\bibitem{kauvo} L.H. Kauffman and P. Vogel, Link polynomials and a
graphical calculus, {\em J. Knot Theory Ramifications\/} {\bf 1} (1992),
59--104.

\bibitem{kawa} A. Kawauchi, {\em A Survey of Knot Theory}, Birkh\"auser,
Boston (1996).

\bibitem{lick} W.B.R. Lickorish, {\em An Introduction to Knot Theory},
Springer, New York (1997).

\bibitem{mena} W. Menasco, Closed incompressible surfaces in alternating
knot and link complements, {\em Topology\/} {\bf 23} (1984), 37--44.

\bibitem{flype} W. Menasco and M. Thistlethwaite, The classification of
alternating links, {\em Ann. of Math.\/} {\bf 138} (1993), 113--171.

\bibitem{mura} K. Murasugi, Jones polynomials and classical conjectures
in knot theory, {\em Topology\/} {\bf 26} (1987), 187--194.

\bibitem{mubook} K. Murasugi, {\em Knot Theory and Its Applications},
Birkh\"auser, Boston (1996).

\bibitem{rol} D. Rolfsen, {\em Knots and Links}, Publish or Perish,
Berkeley, 2nd printing (1990).

\bibitem{four-reg} J. Sawollek, Embeddings of 4-Regular Graphs into
3-Space, {\em J. Knot Theory Ramifications\/} {\bf 6} (1997), 727--749.

\bibitem{alt} J. Sawollek, Alternating Diagrams of 4-Regular Graphs in
3-Space, {\em Topology Appl.} {\bf 93} (1999), 261--273.

\bibitem{schr} A. Schrijver, Tait's flyping conjecture for
well-connected links, {\em J. Combin. Theory Ser. B\/} {\bf 58} (1993),
65--146.

\bibitem{tait} P.G. Tait, On knots I, II, III, {\em Scientific Papers},
Vol. I, Cambridge Univ. Press, London (1898), 273--347.

\bibitem{this} M.B. Thistlethwaite, Kauffman's polynomial and
alternating links, {\em Topology\/} {\bf 27} (1988), 311--318.

\bibitem{yama} S. Yamada, An Invariant of Spatial Graphs, {\em J. Graph
Theory\/} {\bf 13} (1989), 537--551.

\bibitem{yet} D.N. Yetter, Category Theoretic Representations of Knotted
Graphs in $\mathbb S^3$, {\em Adv. Math.\/} {\bf 77} (1989), 137--155.

\bibitem{yoko} Y. Yokota, Topological invariants of graphs in 3-space,
{\em Topology\/} {\bf 35} (1996), 77--87.

\end{thebibliography}
\end{document}